\documentclass[review,10pt]{elsarticle}

\usepackage[usenames]{color}
\usepackage{fancyvrb} 
\usepackage{bm}
\usepackage{amsmath}
\usepackage{subfigure}
\usepackage{multirow}
\usepackage{graphicx}
\usepackage{amssymb}
\usepackage{array}
\usepackage{url}
\usepackage{rotating}
\usepackage{multicol,lipsum}
\usepackage{algpseudocode}
\usepackage{algorithm}
\usepackage{listings}
\definecolor{mygreen}{RGB}{28,172,0} 
\definecolor{mylilas}{RGB}{170,55,241}

\begin{document}

\begin{frontmatter}
\title{A Novel Approach for Optimal Trajectory Design with Multiple Operation Modes of Propulsion System, Part 1}



\author{Ehsan Taheri \fnref{label1f}}  
\fntext[label1f]{Assistant Professor, Email: etaheri@auburn.edu}
\address{Department of Aerospace Engineering, Auburn University\\ Auburn, AL 36849}
\author{John L. Junkins \fnref{label2f}}
\fntext[label2f]{Distinguished Professor, Email: junkins@tamu.edu}
\address{Department of Aerospace Engineering, Texas A\&M University\\  College Station, TX 77843-3141}
\author{Ilya Kolmanovsky \fnref{label3f}} \author{Anouck Girard \fnref{label4f}}
\fntext[label3f]{Professor, Email: ilya@umich.edu}
\fntext[label4f]{Associate Professor, Email: anouck@umich.edu}
\address{Department of Aerospace Engineering, University of Michigan\\  Ann Arbor, MI 48109-2140}

\begin{abstract}
Efficient performance of a number of engineering systems is achieved through different modes of operation - yielding systems described as ``hybrid'', containing both real-valued and discrete decision variables. Prominent examples of such systems, in space applications, could be spacecraft equipped with 1) a variable-$I_{\text{sp}}$, variable-thrust engine or 2) multiple engines each capable of switching on/off independently. To alleviate the challenges that arise when an indirect optimization method is used, a new framework --- Composite Smooth Control (CSC) --- is proposed that seeks smoothness over the entire spectrum of distinct control inputs. A salient aftermath of the application of the CSC framework is that the original multi-point boundary-value problem can be treated as a two-point boundary-value problem with smooth, differentiable control inputs; the latter is notably easier to solve, yet can be made to accurately approximate the former hybrid problem. The utility of the CSC framework is demonstrated through a multi-year, multi-revolution heliocentric fuel-optimal trajectory for a spacecraft equipped with a variable-$I_{\text{sp}}$, variable-thrust engine. 
\end{abstract}

\end{frontmatter}

\section{Introduction} \label{sec:intro}
Trajectory optimization is an important discipline and plays a pivotal role in designing flight paths for atmospheric flight \cite{kamyar2014aircraft, davoudi2019quad}. In space applications, the recent breakthroughs that have taken place in Solar-Electric Propulsion (SEP) systems have led to unprecedented opportunities in deep space missions \cite{rayman2002design} including missions to multiple asteroids \cite{rayman2006dawn}. Electric thrusters operate at a higher level of efficiency compared to chemical rockets, which leads to delivering larger payloads and the ability to accomplish a diverse class of missions \cite{mengali2008optimal,duchek2015solar,genta2016optimal}, specifically, for multiple small-body rendezvous missions. Unsurprisingly, SEP is now envisioned as a reliable option for in-space logistics supply purposes and for cargo missions \cite{jagannatha2018optimization,woolley2019cargo}.

As a consequence of the evolution of propulsion technology and the desire to reduce the cost of ever more complex missions designs, trajectory optimization remains an active field of research. A significant amount of research has been devoted to designing low-thrust space trajectories using direct methods, indirect methods \cite{conway2012survey} or variants of these \cite{yang2007earth}. A variety of tools have been developed that are capable of solving complex interplanetary problems with various low- to medium- to high-fidelity models for propulsion systems and gravitational forces \cite{sauer1997solar,sims1999preliminary,whiffen2006mystic,englander2012automated}. These tools use various mathematical formulations to determine optimal trajectories based on both direct and indirect optimization methods \cite{betts1998survey,trelat2012optimal}. A fairly comprehensive review of the models, objective functions, and solution approaches commonly used for spacecraft trajectory optimization is conducted in \cite{shirazi2018spacecraft}.

Efficient performance of a number of engineering systems with a family of discrete coupled actuation systems is achieved by treating them as hybrid/switching systems with different modes of operation. The focus of this two-part paper series is to alleviate and overcome some challenges that arise when indirect optimization methods are used to characterize optimal trajectories of such systems. Specifically, the task of generating low-thrust fuel-optimal trajectories for spacecraft with multiple modes of operation of propulsion system is considered. 



Two points are worthy of explanation: first, compared to a number of optimal trajectory design problems (e.g., re-entry trajectory optimization  \cite{grant2014rapid}), heliocentric low-thrust trajectory optimization can exploit reasonably accurate and well-behaved models dominated by Keplerian motion. The second important point is that the forces imparted on a spacecraft (either due to SEP system or other types of perturbations) alter its trajectory only after the effects are accumulated over a relatively long period. As a result, the sensitivity of the final states to local variations of the control is usually small and therefore we may anticipate issues with solution's uniqueness and convergence. 

The ultimate goal, of course, is to characterize the time history of control (e.g., the thrust vector) such that the spacecraft eventually reaches its destination in an efficient manner (e.g., to obtain fuel- or time-optimal solutions) while satisfying some operational constraints. In some cases, this may involve the addition of multiple gravity-assist maneuvers and inclusion of inequality state constraints; such problems are highly non-linear and special techniques are needed for solving them \cite{jiang2012practical,pan2018quadratic}. It is possible to approach the solution of these problems by resorting to hybrid optimal control methods \cite{yam2011low,chilan2013automated,rasotto2016multi,taheri2018shaping,englander2016automated}. To focus on key issues addressed in this paper, we do not consider gravity-assist maneuvers or state inequality constraints. 

Despite inherent high sensitivities, it is frequently possible to make the resulting two-point, boundary-value problems (TPBVPs) amenable to numerical solution. This is achieved by constructing a one- or multi-parameter family of neighboring OCPs (i.e., a homotopic map) that contains both the original, difficult problem and an easy-to-solve problem \cite{kim2014optimal}. The easier problem is first formulated and solved. Then, by varying the embedding parameter, a sequence of neighboring TPBVPs is formed in which neighboring converged solutions are used recursively to initiate the next set of neighboring iterations. The family is constructed such that the final problem, in this generated sequence, is the original problem we wish to solve. There are multiple ways to define the sequence of neighboring problems (frequently a one-parameter family of problems), for instance, by 
\begin{enumerate}
    \item identifying the nonlinear part of the dynamics (i.e., the acceleration due to the central body) and by removing (or attenuating) its contribution in order to define a simpler starting problem \cite{thorne1996approximate,petukhov2012method}. Then, the nonlinear effects are introduced to the problem formulation by, for instance, sweeping a multiplicative parameter. For example, a time-fixed fuel-optimal trajectory design from a geostationary orbit to an L1 Halo orbit in the Earth-Moon restricted three-body model is studied in \cite{taheri2018generic}, where the contribution of moon's gravity is introduced in a homoptopic manner,
    \item modifying the boundary conditions from those available solutions, swept to the boundary condition of interest \cite{petukhov2012method,kim2005continuous,hou2012optimization},
    \item modifying the running cost by introducing quadratic/logarithmic terms such that an easier problem is first solved \cite{bertrand2002new,taheri2016enhanced,taheri2018performance},
    \item altering the dynamics by introducing error and error-control terms into the dynamics \cite{mall2017epsilon},
    \item applying a one-parameter smoothing method directly to the admissible control input \cite{petukhov2012method,rasotto2016multi,taheriagenericsmoothing2018}.
\end{enumerate}

While global convergence can rarely be guaranteed when using some subset of the above, several homotopy approaches have been recently established and shown to work well for most of the regular fuel-optimal problems \cite{zhu2017solving}, i.e., those trajectories that do not involve singular control arcs. For instance, a probability-one homotopy method has been proposed for solving the TPBVP associated with time-optimal trajectories \cite{pan2018new}.

However, there are well-known difficulties that may arise when homotopy strategies are adopted to solve optimal control problems (OCPs) \cite{pan2016double}. The simplest homotopy procedures are known to encounter difficulties for problems that consist of many revolutions, especially if a single-shooting scheme is used for solving the resulting TPBVPs. Thus, multiple-shooting schemes have been developed \cite{stoer2013introduction,meng2018low}. A salient aspect of indirect optimization methods is that when they converge accurately, they result in high-resolution solutions, which is quite instrumental for validation and/or training methods that rely on machine learning algorithms to approximate optimal trajectories \cite{sanchez2018real}. Superior speed performance is achievable for some problems when indirect methods are used \cite{uebel2017computationally}.

Despite these difficulties, indirect methods' necessary conditions provide invaluable insights into the structure of the control \textit{prior} to solving the resulting TPBVPs. For instance, the optimal direction of thrusting (impulsive or continuous) for orbit transfer is governed by part of the costate dynamics associated with the velocity vector - also known as the \textit{Primer Vector} due to Lawden \cite{lawden1963optimal}. However, even having such ``structure of the control'' insights, it remains a difficult numerical challenge to find the particular missing costate boundary conditions that ``drive'' the system to the extremal of interest.

A framework has been recently developed that achieves a continuous representation of multi-phase systems with various force models for state dynamics \cite{saranathan2017relaxed}. Motivated in part by the method of \cite{saranathan2017relaxed}, we present here a new framework, called Composite Smooth Control (CSC), within the indirect optimization formalism to deal with the problem of interplanetary trajectory design. Specifically, we focus on problems in which the spacecraft is equipped either with 1) a variable-$I_{\text{sp}}$, variable-thrust (VIVT) engine (in part 1 of this set of two papers) or 2) a cluster of engines (in part 2 \cite{taheri2019anovelpart2}). However, we treat the problem by modifying the admissible control in the optimal control formulation and introduce an approach to smooth the otherwise discrete switches. We make use of a recently introduced approach --- Hyperbolic Tangent Smoothing (HTS) --- to smoothing \cite{taheriagenericsmoothing2018}, which has utility in two respects: 1) it allows us to vary sharpness of the smooth control switches in a homotopic manner, and 2) we can use the HTS to connect different controls and construct a continuously composite control. The switches are still present, but with user control over the sharpness of these switches. The result is a composite continuous control that captures the entire spectrum of admissible controls. 

A prominent feature of the proposed framework is that not only the optimal instances of transition between different operation modes (in the case of a VIVT engine), but also the optimal number of engines (e.g., in the case of a cluster of engines) as well as their optimal operating conditions are revealed without \textit{a priori} assumptions. The multi-engine problem is considered in part 2 \cite{taheri2019anovelpart2}. 

In the proposed framework, bang-bang control profiles of thrusters are also incorporated, which distinguishes the CSC framework from Ref \cite{saranathan2017relaxed}. There is an important implementation subtlety when formulating and solving problems that consist of bang-bang type controls: \textit{the constraint that determines the activation of the bang-bang control is the so-called switching function associated with that specific control}. In other words, for each control input that has a bang-bang control structure, the switching function, $S$, serves as the constraint ($S = 0$) and is considered as the argument of the HTS method. Therefore, the time associated with this constraint is obtained implicitly, and in an autonomous manner. This means that the time of control switches (between the two extreme limits) are also determined autonomously. Collectively, in addition to incorporating smooth transitions due to multiple time- or state-triggered constraints, CSC enables us to incorporate switching-function-triggered constraints that govern bang-bang type control inputs.

The smoothness of the mathematically continuous composite control input is found to significantly enhance the convergence performance of numerical algorithms that rely on adaptive numerical differential solvers and Newton-type update schemes, which are most commonly used in solving TPBVPs. Using all of the above considerations, we are able to formulate and solve low-thrust trajectory optimization problems involving interplanetary heliocentric phases of flight when a single-shooting solution scheme is used and with enhanced reliability and efficiency. In fact, we use MALTAB's built-in numerical solver, \textit{fsolve}, and numerical integration function, \textit{ode45}, for all of our numerical simulations to make our results readily reproducible. 

The CSC framework allows us to solve such problems in 1) a unified manner, 2) without resorting to different algorithms, and 3) avoid switch time determination iterative procedures during dynamical motion integration, which is an essential step in traditional methods to accurately solve such multi-mode systems. 

A second contribution of the paper is an alternative approach to model the contribution of the planetary gravitational perturbations in the costate differential equations when the set of modified equinoctial elements (MEEs) are used to represent the trajectory of the spacecraft. This is achieved through a complex-based numerical derivative approach to determine the Jacobian of the Hamiltonian that is fundamental to deriving the dynamics of the costates. Key implementation details of this alternative, effective, numerical approach are highlighted.  

The remainder of the paper is organized as follows. Section \ref{sec:models} reviews the equations of motions and summarizes some of the important factors involved in modeling a VIVT engine as well as solar array and power sub-systems. Section \ref{sec:variableIspvariablethrust} presents the OCP for a spacecraft that is equipped with a VIVT engine. Implementation details of the complex-based derivative method are explained in Section \ref{sec:variableIspvariablethrust}. Section \ref{sec:CSC} presents the CSC framework and its application to solving fuel-optimal trajectory optimization problems. Section \ref{sec:results} presents the results. Finally, Section \ref{sec:conclusion} provides a conclusion to the paper.  

\section{Power System, Thruster and Perturbation Modeling} \label{sec:models}
In this section, we review the solar electric power models and the disturbing accelerations due to secondary bodies. Equations of motion are given while taking into account the performance of power sub-system and actuation models.  
\subsection{Power Sub-system Modeling}
Application of SEP systems for planetary missions is not a new concept \cite{sauer1987application}. A crucial step in high-fidelity trajectory design for spacecraft that rely on SEP propulsion is to develop and use appropriate power models since it very significantly affects the performance of the engines and the ensuing trajectories. Solar arrays are responsible for powering the spacecraft and sustaining the operation of its various sub-systems during the entire time of flight. The SEP system, of course, may or may not be providing thrust at any point in time along an optimal trajectory. In some cases, there are multiple engines that can be on (or not) and in others, the engine design allows for variable $I_{\text{sp}}$, throttle thrust, and propellant consumption. 

Depending on their size, solar arrays generate a nominal beginning-of-life power, $P_{0,\text{BOL}}$, typically defined at launch time as the power at one astronomical unit (AU) from the Sun. For a fixed-size solar array, the nominal power value can change due to two reasons: 1) the distance of the spacecraft from the Sun changes, and 2) the efficiency of the solar cells degrades over time due to thermal cycling, plume impingement, radiation, and other environmental effects. 

Let $\textbf{r}$ denote the position vector of the spacecraft expressed in any Sun-centered inertial frame of reference, and let $r = \| \textbf{r}\|$ denote the Euclidean norm of the position vector, that is, the distance to the Sun. Then, the fractional multiplicative variation of the nominal power due to the aforementioned sets of physical effects can be approximated \cite{ellison2018application} by introducing distance- and time-dependent multiplicative functions as
\begin{align} 
\phi(r) = &   \frac{1}{r^2} \left [ \frac{A_1+\frac{A_2}{r}+\frac{A_3}{r^2}}{1+A_4 r +A_5 r^2} \right ], \label{eq:radterm}\\
\psi(t) = & (1-\sigma)^{\tau} \label{eq:timeterm}, 
\end{align}
where $\phi(r)$ denotes the distance-dependent term, $\psi(t)$ denotes the time-dependent term (due to aging), $\sigma$ is decay rate of the solar arrays measured as a percentage per year (usually between 2 to 4 \%/year) and $\tau = \tau(t)$ denotes elapsed time from the launch time measured in years \cite{ellison2018application}. The bracketed term in Eq.~\eqref{eq:radterm} is a quasi-empirical term, which, with approximate solar power system-specific choices of the five $A_i$ coefficients, approximates the solar intensity induced variations of the output power from the array. Requiring $\phi(r) = 1$ at $r = 1$ AU gives rise to the linear constraint equation 
\begin{equation*}
    A_1 + A_2 +A_3 -A_4 -A_5 = 1,
\end{equation*}
which can be imposed on the fitting processor, thereby eliminating one of $A$ coefficients in Eq.~\eqref{eq:radterm} and reduce the unknowns from 5 to 4. While this approximation may differ typically from the truth by a fraction of one percent, we adopt Eqs.~\eqref{eq:radterm} and \eqref{eq:timeterm} to capture representative radial and time variability of the solar power available with the coefficients $A_i$ ($i = 1, \cdots, 5$), which are determined through fitting measurements from ground tests in laboratories. Ultimately, the power generated by solar arrays is approximated herein as
\begin{equation}
P_{\text{SA}}(t,r) = \psi(t) \phi(r) P_{0,\text{BOL}},
\end{equation}
where $P_{0,\text{BOL}}$ corresponds to power at $t = t_0$. The power needed to sustain the operation of sub-systems of a spacecraft as well as the power needed for operation of Power Processing Unit (PPU) can also be approximated \cite{ellison2018application} as a reciprocal power function of $r$
\begin{equation} \label{eq:psc}
P_{s/c} = D_1 + \frac{D_2}{r} +\frac{D_3}{r^2}.
\end{equation}

\begin{figure}[htbp!]
\centering
\includegraphics[width=3.0in]{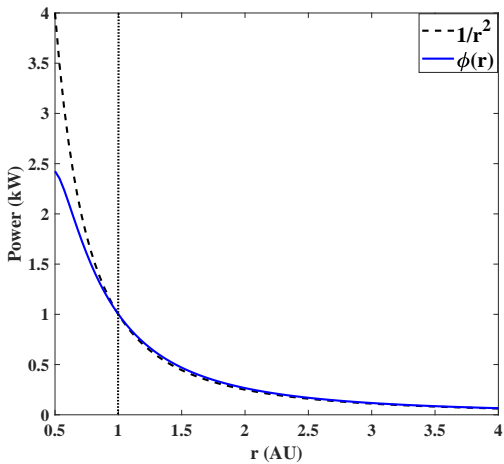}
\caption{Variation of $P_{\text{SA}}$ vs. distance to the Sun for $P_{0,\text{BOL}} = 1$ kW; $\psi(t) = 1$.}
\label{fig:powervsr}
\end{figure}

The PPU modulates the power sent to each thruster. The coefficients in Eq.~\eqref{eq:psc} are difficult to determine since power consumption is a complicated function of not only the heat transfer, which depends on the distance to the Sun and available solar power intensity, but also on the orientation of the spacecraft. Only the dependency on distance is modeled in Eq.~\eqref{eq:psc}. For numerical results in this paper, $P_{s/c} = P_{\text{ppu}}$ and $P_{\text{ppu}}$ is assumed to be the required power to operate all spacecraft sub-systems other than the thruster(s). The coefficients $D_j$ ($j = 1,\cdots,3$) are determined through simulated tests in ground-based laboratories. If a sufficiently precise on-board acceleration can be measured, in principle, a parametric representation analogous to Eqs.~\eqref{eq:radterm} and \eqref{eq:timeterm} could be updated during the mission.  The available power to be used for thrusters then becomes
\begin{equation} \label{eq:avapower}
P_{\text{ava}} = P_{\text{SA}}(t,r) - P_{s/c}.
\end{equation}
Figure \ref{fig:powervsr} depicts the variation of solar array power, $P_{\text{SA}}$, versus distance to the Sun when 1) only an inverse squared dependency to $r$, i.e., $\phi(r) = 1/r^2$ is considered, and 2) the following values are considered for the coefficients of $\phi(r)$: $A_1 = 1.321$, $A_2 = -0.108$, $A_3 = -0.117$, $A_4 = 0.108$, $A_5 = -0.013$ \cite{laipert2015automated}. 
It is clear that at distances greater than 1 AU (to the right of the vertical dotted line) the differences between the two approximations are negligible; however, for trajectories to inner planets (to the left of the vertical dotted line) or for multi-revolution transfer orbit that passes inside 1 AU, more accurate models have to be used, otherwise errors larger than 25\% may occur \cite{sauer1997solar}. 
 
\subsection{Variable-$I_{\text{sp}}$, Variable-Thrust Engine Model}
In this work, we consider a specific type of engine that is capable of modifying its specific impulse, $I_{\text{sp}}$, and thrust level, $T$, simultaneously. The use of VIVT engines is shown to impact trip times and propellant requirements \cite{chang1995rapid,kechichian1995optimal,seywald2003fuel,ranieri2005optimization}. This type of engine is considered so that we can demonstrate the utility of the proposed indirect optimization framework (since the optimal values of specific impulse can be expected to switch between different admissible values \cite{casalino2004optimization}). 

The input power and the exhaust velocity are adapted as control variables. For a VIVT engine, thrust magnitude, $T$, is expressed as 
\begin{equation} \label{eq:thrustVIVT}
T = \frac{2 \eta P}{c},
\end{equation}
where $\eta$ is the assumed constant engine efficiency, $P$ is the power input, and $c = I_{\text{sp}} g_0$ is the exhaust velocity, $I_{\text{sp}}$ is the specific impulse, and $g_0$ is the Earth's gravitational acceleration at sea level. Let $\hat{\bm{\alpha}}$ denotes the thrust steering unit vector, and let $m$ denote the mass of spacecraft; the propulsion acceleration vector can be written as 
\begin{equation} \label{eq:VIVT_u}
\textbf{u}_{\text{prop}} = \frac{T}{m} \hat{\bm{\alpha}} = \frac{1}{m}\frac{2 \eta P}{c} \hat{\bm{\alpha}}.
\end{equation}

This particular parameterization of the total thrust acceleration vector (i.e., a unit vector, $\hat{\bm{\alpha}}$, and the magnitude of thrust, $T$) has its advantages, which lends itself properly to Lawden's primer vector theory. As a consequence of the thrust expression given in Eq.~\eqref{eq:thrustVIVT}, the time rate of change of mass can be written as a function of $(P,c)$ as
\begin{equation} \label{eq:VIVT_mdot}
\dot{m} = -\frac{T}{c} = -\frac{2 \eta P}{c^2}.
\end{equation}

Expressions given in Eqs.~\eqref{eq:VIVT_u} and \eqref{eq:VIVT_mdot} will be used in deriving the equations of motion. When adopting ($P,c$) as control variables in the context of any particular engine, they can be considered real variables that must lie in bounded sets
\begin{align} \label{eq:controlbounds}
\{ P_{\text{min}} \leq P \leq P_{\text{max}} \};~~~~~\{ c_{\text{min}} \leq c \leq c_{\text{max}}\}. %
\end{align}

\subsection{Secondary-Body Perturbation Modeling}
For heliocentric phase of SEP-based flight, a two-body model is considered accurate enough for preliminary analysis. For instance, it is explained that for Deep Space 1 mission \cite{rayman2002design}, the preliminary optimal solutions obtained from Solar Electric Propulsion Trajectory Optimization Program (SEPTOP) were further optimized through a secondary tool, NAVTRAJ, that use a higher-fidelity gravitational perturbation model, including secondary bodies and solar radiation perturbations. 

According to \cite{rayman2002design}, for some cases, the errors due to inaccurate modeling of the power system had a much greater offset on the solution in comparison to ignoring higher-fidelity gravitational perturbations. Nevertheless, we are interested in quantifying the importance of these perturbations on the considered test cases in the present study. In addition, the details of deriving the costate dynamics are presented for two approaches when secondary-body perturbations are incorporated into the state dynamics. 

Let $\textbf{r}$ and $\textbf{v}$ denote the position and velocity vectors of the spacecraft, respectively, in an inertial frame of reference ($I$). The acceleration due to secondary bodies, $\textbf{a}_{\text{sb}}$, can be evaluated as \cite{betts2015optimal}
\begin{equation} \label{eq:secondaryacc}
\textbf{a}_{\text{sb}} = - \textbf{C}_{I}^{\text{LVLH}} \sum_{j=1}^{N_{\text{sb}}} \mu_{j} \left [ \frac{ \textbf{r} - \textbf{r}_{\odot}}{\| \textbf{r} - \textbf{r}_{\odot} \| ^3 } + \frac{ \textbf{r} - \textbf{r}_j}{\| \textbf{r} - \textbf{r}_j \| ^3 } \right ],
\end{equation}
where $N_{\text{sb}}$ is the number of considered secondary bodies, $\textbf{C}_{I}^{\text{LVLH}}$ is the orthogonal transformation matrix that projects the acceleration vector components in the inertial frame of reference to acceleration vector components in the Local Vertical Local Horizontal (LVLH) frame and $\mu_j$ denotes the gravitational parameter of the $j$-th secondary body. Let $\hat{\textbf{u}}_r = \textbf{r}/\| \textbf{r}\|$, $\hat{\textbf{u}}_h = \textbf{r} \times \textbf{v} /\| \textbf{r} \times \textbf{v} \|$, and $\hat{\textbf{u}}_t =  \hat{\textbf{u}}_h \times \hat{\textbf{u}}_r$ be the expressions for the osculating unit vectors of the LVLH frame. At each instant of time, the transformation (direction cosine) matrix can be represented as the three orthogonal LVLH unit vectors as  
\begin{equation}
\textbf{C}_{I}^{\text{LVLH}} = \left [ \hat{\textbf{u}}_r, \hat{\textbf{u}}_t, \hat{\textbf{u}}_h \right ]^{\top}.
\end{equation}

\subsection{Equations of Motion}
In the heliocentric phase of flight, the spacecraft motion is predominantly governed by the gravitational attraction of the Sun. In addition, the spacecraft is equipped with an electric thruster or a cluster of thrusters and its mass, $m$, changes due to the consumption of propellant. We choose to formulate an OCP in terms of the modified equinoctial orbital elements (MEEs) \cite{walker1986set} since it has been shown that these elements are very attractive for describing low-thrust trajectories \cite{taheri2016enhanced,junkins2018exploration}. We have separated the dynamics of MEEs and time rate of change of spacecraft mass to facilitate derivation of optimality conditions. Let $\textbf{x} = [p,f,g,h,k,l]^{\top}$ denote the vector of MEEs and let $\textbf{a} = [a_r, a_t, a_n]^T$ represent the non-two-body accelerations expressed in the LVLH frame acting upon the spacecraft. The dynamics of MEEs can be written as
\begin{align} 
\dot{\textbf{x}} =& \textbf{A}(\textbf{x},t) + \mathbb{B}(\textbf{x},t) \textbf{a},
\end{align}
where $\textbf{A} \in \mathbb{R}^{6 \times 1}$ denotes the unforced vector part of the MEEs dynamics and $\mathbb{B} \in \mathbb{R}^{6 \times 3}$ is the control influence matrix defined in \cite{junkins2018exploration}.

The MEE set has five slow variables and one (very regular) fast variable, $l$. The complete dynamical model including the variation of mass becomes 
\begin{align}
\dot{\textbf{x}} = & \textbf{A}(\textbf{x},t) + \mathbb{B}(\textbf{x},t) \textbf{a}, \label{eq:elementsdiffeq}\\ 
\dot m = & -\frac{2 \eta P}{c^2}. \label{eq:massdiffeq}
\end{align}

In Eq.~\eqref{eq:elementsdiffeq}, the acceleration vector, $\textbf{a}$, can be conveniently expressed as the sum of two acceleration terms as 
\begin{equation}
\textbf{a} = \textbf{u}_{\text{prop}} + \textbf{a}_{\text{sb}},
\end{equation}
where $\textbf{u}_{\text{prop}}$ denotes the control acceleration vector due to engine(s) (see Eq.~\eqref{eq:VIVT_u}) and $\textbf{a}_{\text{sb}}$ denotes the perturbing acceleration vector due to all secondary bodies (i.e., the other planets in the Solar System) derived in Eq.~\eqref{eq:secondaryacc}. There exist other smaller disturbances such as solar radiation pressure that are not considered in this work. 

We mention that we have previously investigated eight choices of coordinates and elements in \cite{junkins2018exploration} for the class of time-fixed rendezvous-type fuel-optimal problems that involve indirect trajectory optimization. The MEEs were shown superior computationally to the seven competing coordinate choices considered \cite{junkins2018exploration}. The relative merits of the coordinate choices, however, may vary due to the magnitude of the controls, the gravitational forces as well as the particular boundary conditions imposed. Furthermore, the geometric insights using some unique choices of coordinates may, in some cases, lead to other decisions on ``which coordinates are the best''.

\section{Fuel-Optimal Indirect Formulation for VIVT Engines} \label{sec:variableIspvariablethrust}
In this section, the TPBVP associated with the fuel-optimal trajectories is developed. To keep this paper reasonably self-contained, most relevant equations are derived. In particular, the Hamiltonian, $H$, plays a central role in optimal control necessary conditions and we refer to the various terms in $H$ to derive the Jacobian needed in the costate dynamics. We explain the complex-based derivative method as an alternative to deriving the costate dynamics.
\subsection{Formulation of the Fuel-Optimal Problem}
For a fuel-optimal problem, the minimization of the cost functional  
\begin{eqnarray} \label{eq:costfunction}
J = -m(t_f),
\end{eqnarray}
is of interest. We form the Hamiltonian associated with the OCP
\begin{equation} \label{eq:VIVTHamil}
H = \bm{\lambda}^{\top} \left [ \textbf{A}(\textbf{x},t) + \frac{1}{m}\frac{2 \eta P}{c} \mathbb{B}(\textbf{x},t)  \hat{\bm{\alpha}} \right ] - \lambda_m \frac{2 \eta P}{c^2},
\end{equation}
where $\bm{\lambda} = [\lambda_p,\lambda_f,\lambda_g,\lambda_h,\lambda_k,\lambda_l]^{\top}$ denote the costate vector associated with the MEEs, $\textbf{x}$ and $\lambda_m$ denotes the costate associated with the mass state, $m$. 

The control inputs ($\hat{\bm{\alpha}}$, $P$, and $c$) have to be characterized such that $H$ is minimized according to the PMP. The first two control inputs (i.e., $\hat{\bm{\alpha}}$ and $P$) appear bi-linearly in $H$, thereby requiring PMP to be used, whereas the optimal value of specific impulse can be obtained through the strong form of PMP (i.e., $\partial H/\partial c = 0$), if the minimizing value of $c$ lies in the admissible range: $c_{\text{min}} < c < c_{\text{max}}$, otherwise, PMP must be invoked and $c$ will lie on the boundary of the admissible region that minimizes $H$. $\|\hat{\bm{\alpha}} \| =1$, but $\hat{\bm{\alpha}}(t)$ is otherwise an arbitrary real vector and all $\{ P(t), c(t) \}$ must lie in a bounded feasible set of Eq.~\eqref{eq:controlbounds}.

In the remainder of the paper, the arguments of $\textbf{A}$ and $\mathbb{B}$ are suppressed to have compact expressions. As is evident by inspection of $H$, the optimal (denoted by superscript `*') direction of thrust that minimizes $H$ over all feasible $\hat{\bm{\alpha}}(t)$ is the unit primer vector
\begin{equation} \label{eq:primervector}
\hat{\bm{\alpha}}^* = -\frac{\mathbb{B}^{\top} \bm{\lambda}}{ \| \mathbb{B}^{\top} \bm{\lambda} \| }.
\end{equation}
Upon substituting the expression for the optimal thrust direction given by Eq.~\eqref{eq:primervector} into Eq.~\eqref{eq:VIVTHamil} and collecting terms, $H$ can be written as
\begin{equation}
H = H_0 -  \left [ \frac{ \| \mathbb{B}^{\top} \bm{\lambda} \|}{m} + \frac{\lambda_m}{c}\right ] \frac{2 \eta P}{c} \Longleftrightarrow H = H_0 -  \left [ \frac{ \| \mathbb{B}^{\top} \bm{\lambda} \|}{m} + \frac{\lambda_m}{c}\right ] T,
\end{equation}
where $H_0$ represents control-independent terms. Note that $H$ is written in two forms: one where power is the control input and one where thrust, $T = 2 \eta/c$, is the control input. The optimal power and/or thrust are obtained as 
\begin{align} \label{eq:VIVTpowercondition}
P^{*}(S) = & \begin{cases}
     P_{\text{max}} & \text{if}~S>0,  \\
     0 & \text{if}~S<0,
\end{cases} & \text{or} & & T^{*}(S) = & \begin{cases}
     T_{\text{max}} & \text{if}~S>0,  \\
     0 & \text{if}~S<0,
\end{cases}
\end{align}
where the power switching function, $S$ is defined as
\begin{equation}\label{eq:sfm-p1}
S \equiv \frac{\|\mathbb{B}^{\top} \pmb{\lambda} \|}{m}+\frac{\lambda_m}{c}.
\end{equation}

We have assumed that singular arcs (where $S = 0$ for a finite time interval) do not occur during the trajectory. The first term of Eq.~\eqref{eq:sfm-p1} is always positive, so $S$ can change sign only when the second term is a sufficiently large negative number relative to the first term, and $S = 0$ switches occur from exact cancellations. For a VIVT engine, $c \in [c_{\text{min}}, c_{\text{max}}]$, where $c = c_{\text{max}}$ leads to the least value of the second term, which in turn results in the maximum value of $S$. So, $c_{\text{max}}$ dictates $P = P_{\text{max}}$ and $T = T_{\text{max}}$. In the numerical results we observed that $\lambda_m(t) <0$ along all extremal solutions. 



We still have to characterize the optimal value of the specific impulse. $H$ is a quadratic function of $1/c$; thus, it is straightforward to apply the strong form of optimality, i.e., $\partial H/\partial c = 0$ and solve for the optimal (minimizing) value of $c$ as
\begin{equation} \label{eq:VIVTcop}
c_{\text{op}} = -\frac{2m \lambda_m}{\|\mathbb{B}^{\top} \pmb{\lambda} \|}.
\end{equation}

Altogether, $H$ is minimized if $c_{\text{op}}$ in Eq.~\eqref{eq:VIVTcop} is used as long as its value lies within the admissible range ($c_{\text{op}} \in [c_{\text{min}}, c_{\text{max}}]$). Clearly, the nearest violated limit is taken if $c_{\text{op}}$ lies outside its admissible range. 

The optimal $c^*$ is summarized as
\begin{align} \label{eq:VIVTcthree}
c^* = \begin{cases}
c_{\text{max}} & \text{if}~ c_{\text{op}} \geq c_{\text{max}},  \\
c_{\text{op}}  & \text{if}~ c_{\text{op}} \in (c_{\text{min}}, c_{\text{max}}), \\
c_{\text{min}} & \text{if}~ c_{\text{op}} \leq c_{\text{min}}.
\end{cases}
\end{align}

The above set of conditions is usually implemented through a \textit{max-min} structure as
\begin{equation} \label{eq:minmax}
c^* = \text{max} \left \{ c_{\text{min}},\text{min}\left [ c_{\text{op}}, c_{\text{max}} \right ] \right \},
\end{equation}
where $\text{max}[.,.]$ ($\text{min}[.,.]$) expression returns the maximum (minimum) of the two arguments. 

Unless the optimal $c$ always satisfies $c_{\text{op}} \in (c_{\text{min}}, c_{\text{max}})$, the structure defined in Eq.~\eqref{eq:minmax} is inherently nonsmooth. This non-smoothness is essentially due to the fact that, at the time of implementation, a set of ``if $\rightarrow$ then'' conditional statements are required to select the optimal values of exhaust velocity. The particular implementation in Eq.~\eqref{eq:minmax} is used in a large class of problems and an approach to smooth this structure (as is presented herein) has a great application for a vast class of problems. 

There is another source of non-smoothness introduced to the problem through optimality criterion that specifies when to switch the thruster ``on'' or ``off'' (see Eq.~\eqref{eq:VIVTpowercondition}). Our goal is to avoid all non-smoothness by constructing smooth transitions between not only the three cases of $c^*$, but the transitions between possible power inputs. The reason we seek smooth approximations of on-off conditions is to make numerical solutions more convenient, efficient and accurate, and as generalized later in this paper to achieve smooth invariant embedding of multi-mode thrust controls. 

We proceed by forming the remaining components of the necessary conditions, namely, the costate dynamics and boundary conditions. At this stage, let $\textbf{a}_{\text{sb}} = \textbf{0}$ for the sake of simplicity of deriving the costate dynamics ($\dot{\bm{\lambda}} = \dot{\bm{\lambda}}_{\text{two-body}}$). From the expression for $H$ given in Eq.~\eqref{eq:VIVTHamil} together with the Euler-Lagrange equation, the dynamics of the costates can be derived as
\begin{align}
\dot{\bm{\lambda}}_{\text{two-body}} = & -\left [ \frac{\partial H}{\partial \textbf{x}} \right ]^{\top} , \label{eq:costatex} \\ 
\dot{\lambda}_m = & -\frac{\partial H}{\partial m} = - \frac{\|\mathbb{B}^{\top} \pmb{\lambda} \|}{m^2} \frac{2 \eta P}{c}. \label{eq:costatemass}
\end{align}

The costate dynamics associated with the MEEs are given in other references  \cite{quarta2011minimum}. However, there is a subtlety involved where $H$ is indirectly a function of the solar distance since the solar array output power varies as a function of distance from the Sun and this dependence has to be taken into consideration when deriving the costate dynamics (see under ``TRAJECTORY OPTIMIZATION'' Section in \cite{sauer1997solar}). 


The mathematical expressions for the right-hand side of the $\dot{\bm{\lambda}}$ equations become lengthier when $\textbf{u}_{\text{sb}} \neq \textbf{0}$. Section \ref{eq:costatederivation} is devoted to present two alternatives for deriving the costate dynamics when perturbations due to secondary bodies are taken into consideration. Note that the inclusion of perturbations due to secondary bodies does not alter the \textit{mathematical expressions} derived for optimal control inputs (namely, Eqs.~\eqref{eq:primervector} through \eqref{eq:minmax}) since these disturbing accelerations are not functions of control and are only functions of position and velocity vectors. Nevertheless, they eventually alter the time history of costates, and indirectly alter the resulting optimal profile of control inputs.
\subsection{Two-Point Boundary-Value Problem}
Depending on the type of maneuver, different boundary conditions can be enforced. In this paper, we are dealing with time-fixed rendezvous-type problems, where the final mass is free. The position and velocity vectors of the spacecraft in rendezvous maneuvers match their target body counterparts (denoted by subscript `T'). It is assumed that the spacecraft leaves the Earth's Sphere of Influence (SOI) on a parabolic trajectory, i.e., with zero excess velocity $v_{\infty} = 0$ and that, on a solar scale, the SOI is negligibly small compared to 1 AU and the inertial velocity of the spacecraft is the inertial velocity of the Earth. This is a frequently used first approximation in the solar orbit transfer preliminary mission design. The final boundary conditions (seven equality constraints) can be written in the vector function form as
\begin{equation}\label{eq:posvelconp1}
\bm{\psi}(\textbf{x}(t_f),t_f) = \left [ [\textbf{x}(t_f) - \textbf{x}_T]^{\top},\lambda_m(t_f) + 1 \right ]^{\top} = \textbf{0},
\end{equation}
where $\textbf{x}_T = [p_T,f_T,g_T,h_T,k_T,l_T]^{\top}$ denotes the final target state values. Since final mass is free, the final value of the mass costate has to be -1 (due to transversality condition). Another unknown of the problem, for multi-revolution trajectories, is the number of en-route revolutions, $N_{\text{rev}}$, in the transfer orbit, which has to be determined. Its value is taken into consideration when the change in the true longitude, $l$, is to be enforced as a boundary condition,
\begin{equation}
l_T = \tilde{l}_f + 2 \pi N_{\text{rev}},
\end{equation}
where $\tilde{l}_f$ is the initial true longitude of the target body (i.e., for $N_{\text{rev}} = 0$). $l_T$ is essentially the modified true longitude boundary condition based on the value of $N_{\text{rev}}$.

\textbf{Remark}: Previous studies have shown that for each specified number of revolutions belonging to the feasible set of $N_{\text{rev}}$ integers, there is one local extremal for \textit{time-fixed rendezvous-type fuel-optimal} trajectories \cite{taheri2018unified}. By seeking extremals for all feasible $N_{\text{rev}}$, we can evidently obtain the global extremal. 

Let $\textbf{z} = [\textbf{x}^{\top},m,\bm{\lambda}^{\top},\lambda_m]^{\top}$ denote the state-costate vector, then, we can write,
\begin{equation} \label{eq:F}
\dot{\textbf{z}} = \textbf{F} = [\dot{\textbf{x}}^{\top}, \dot{m}, \dot{\bm{\lambda}}^{\top}_{\text{two-body}}, \dot{\lambda}_m]^{\top},
\end{equation}
where $\bm{\alpha} = \bm{\alpha}^*$, $P = P^*(S)$ and $c = c^*(S)$ (Note $\dot{\textbf{x}}$, $\dot{m}$, $\dot{\bm{\lambda}}_{\text{two-body}}$, $\dot{\lambda}_m$ here are shorthand for the RHS of Eqs.~\eqref{eq:elementsdiffeq},\eqref{eq:massdiffeq},\eqref{eq:costatex},\eqref{eq:costatemass}). Once the optimal values of the control components are substituted into $\textbf{F}$, the equations of motion can be integrated numerically, if initial conditions are fully specified. 

In the class of orbit maneuvers we consider, only the initial state $\textbf{x}(t_0) = \textbf{x}_0$ and $m(t_0) = m_0$ are specified. The final state $\textbf{x}(t_f)$ as well as the final costates are functions of the initial costate $\bm{\eta}(t_0)$, where $\bm{\eta}(t_0) = [\bm{\lambda}^{\top}(t_0),\lambda_m(t_0)]^{\top}$ is the vector of unknown variables to be determined such that Eq.~\eqref{eq:posvelconp1} is satisfied. Thus, we have a TPBVP that can be solved if we have a starting estimate $\bm{\eta}(t_0)$ within the domain of convergence of the particular algorithm to satisfy the prescribed boundary conditions. There are seven constraints in Eq.~\eqref{eq:posvelconp1} and seven unknown elements in $\bm{\eta}(t_0)$. In the next section, the details of applying the composite smoothing method for solving the resulting TPBVP is presented. We mention that, the much slower time rate of change and near-linearity of the MEEs compensate for the more algebraically complicated expressions, as compared to traditional Cartesian or spherical coordinates \cite{taheri2016enhanced,junkins2018exploration}. 

\subsection{Derivation of Costate Dynamics} \label{eq:costatederivation}
In indirect optimization methods, costates constitute a crucial part of the necessary information to construct the optimal control. However, derivation and coding of the governing costate dynamics can become quite tedious and is an error-prone process for complicated systems such as the one in this study. Moreover, the choice of coordinates or elements usually has a significant impact on both the derivation and the solution process when the low-thrust trajectory optimization problem is considered \cite{taheri2016enhanced,junkins2018exploration}. The simplest expressions for differential equations of the costates arise when the set of Cartesian coordinates is used to represent state dynamics. The expressions of the  differential equations of the costates become lengthier and more involved when the set of MEEs is utilized. The relations become even more involved with the inclusion of gravitational perturbations such as the zonal harmonics term, $J_2$ \cite{kechichian2007streamlined}, and especially so for higher order series representations of gravitational perturbations. We mention in passing that the complexity of deriving the costate dynamics (for a particular coordinate choice) does not generally correlate to the efficiency and reliability of numerical convergence \cite{junkins2018exploration}.

\subsubsection{Symbolic Approach}
The error-prone process of deriving the costate differential equations can be avoided by using symbolic toolboxes such as MATLAB' Symbolic Math Toolbox to construct the necessary expressions \cite{taheri2016enhanced,taheri2018performance} and to a significant degree, automate code generation. For example, an automated algorithm can be devised for generating the expressions as well as generating functions (using MATLAB's \textit{matlabFunction} command that converts symbolic expressions to a function handle or a file) that can be readily used for propagating the costate dynamics \cite{taheri2016enhanced}.

Although the algebraic process becomes automated, the generation of files may take a significant computational time and the functional expressions may not always be as well optimized for subsequent numerical computation as equivalent expressions derived and programmed by hand. A portion of the time delay is caused by the fact that \textit{matlabFunction} looks for repeating expressions and factorization within the symbolic expressions so that intermediate variables are introduced to avoid repeated evaluations of such expressions and to improve the computational efficiency. However, when perturbing accelerations are included in state dynamics, the generation of a reasonably efficient code can take significant time. Even when the generation of symbolic files is successful, the resulting expressions may ``explode'' and become so lengthy that evaluation of such expressions can negate computational advantages of using symbolic expressions. There is an alternative \textit{numerical} approach for determining the costate dynamics to near-machine precision, which will be discussed in Section \ref{sec:complexapproach}.

For the considered problem, the main disturbances are due to secondary bodies and the contribution of these disturbances to the Hamiltonian, $H_{\textbf{sb}}$, can be written as
\begin{equation}
H_{\textbf{sb}} = \bm{\lambda}^{\top} \mathbb{B}(\textbf{x},t)  \textbf{a}_{\text{sb}}(\textbf{x},t),
\end{equation}
where $\textbf{a}_{\text{sb}}$ is defined in Eq.~\eqref{eq:secondaryacc}. Therefore, it is possible to treat the derivation of costate dynamics in two steps: 1) the gravitational acceleration contribution due to the Sun and control input, which is formulated in the previous section, and 2) the gravitational acceleration contribution due to the secondary bodies. The contribution of the $\upsilon$-th secondary body to the costate dynamics can be derived as
\begin{equation}
\dot{\bm{\lambda}}_{\text{sb},\upsilon} =  - \left [ \frac{H_{\textbf{sb}}}{\partial \textbf{x}} \right ]^{\top} = - \sum_{j=1}^3 \left \{ \left [ \frac{\partial \tilde{\mathbb{B}}_j}{\partial \textbf{x}} \right ]^{\top} \bm{\lambda} ~\textbf{a}_{\text{sb},j} \right \}  - \left [ \frac{\partial \textbf{a}_{\text{sb}}}{\partial \textbf{x}} \right ]^{\top} \mathbb{B}^{\top} \bm{\lambda},
\end{equation}
where $\textbf{}$ $\tilde{\mathbb{B}}_j$ denotes the $j$-th row of $\textbf{a}_{\text{sb}}$ in Eq.~\eqref{eq:secondaryacc}. The symbolic expressions can be generated once and used for other planetary bodies. 

The explicit expressions for $\partial \textbf{a}_{\text{sb}}/\partial \textbf{x}$ are rather involved since they depend on the transformation from the inertial to the LVLH frame, which, in turn, depends on the position and velocity vectors as well as the perturbation term, all of which have a non-linear dependence on $\textbf{x}$. There are certain aspects of the problem that can be exploited to develop faster algorithms. For instance, the numerical value of the transformation matrix remains the same at each time instant and irrespective of the number of planetary bodies; hence, at each time instant, its value can be calculated and used whenever it is required. Consequently, the total contribution due to secondary bodies can be written as
\begin{equation}
\dot{\bm{\lambda}}_{\text{sb}} = \sum_{i = 1}^{N_{\textbf{sb}}} \dot{\bm{\lambda}}_{\text{sb},i}.
\end{equation}

The contribution of secondary bodies are added to the central body (i.e., two-body gravitational model), where the costate dynamics associated with the MEEs become
\begin{equation}
\dot{\bm{\lambda}} = \dot{\bm{\lambda}}_{\text{two-body}} + \dot{\bm{\lambda}}_{\text{sb}}.
\end{equation}

Secondary bodies do not contribute to the costate associated with the mass since we deal with restricted dynamics, i.e., ignoring the mass of the spacecraft compared to that of planets.

\subsubsection{Numerical Approach Using Complex-Based Derivative Approach} \label{sec:complexapproach}
There is a second approach to generate the time history of costates through the complex-based derivative (CX) method \cite{martins2013review}. The application of this method for generating the State Transition Matrix to be used for facilitating the solution of fuel-optimal low-thrust trajectories is explained in \cite{taheri2016enhanced} (see Section III in \cite{taheri2016enhanced}). 

While being computationally more demanding (typically) compared to the numerical computation with the symbolic expressions, this approach presents distinct advantages with respect to the elapsed time for human effort in the implementation process. The code is also much more compact. In some problems, such a numerical method becomes faster than the evaluation of lengthy symbolic expressions of the costates dynamics. This is especially true if one uses MATLAB, which is an interpretative programming environment. In addition, it is difficult to apply the symbolic approach to problems in which data tables are used for interpolation purposes. However, the CX method handles such cases as a black-box and can usually accommodate interpolation of data tables, although occasional anomalies arise if iterations cause solutions to converge near table entries where derivatives may be discontinuous (depending on the interpolation process utilized). These problems can be avoided by using Lagrange interpolation centered on the nearest table entry and considering more local tabulated data points for support. In addition, it is possible to have complicated dynamics if one is interested to investigate the impact of higher-order gravity models for low-altitude planetocentric trajectory design.

Here, the procedure is explained for obtaining the numerical Jacobian value for one representative costate. The same procedure can be followed for other costates of the system. Let $\lambda$ be one of the costates of the problem. Our goal is to evaluate the \textit{numerical value} of the time rate of change of the respective costate, $\dot{\lambda}$, which can be obtained as
\begin{equation} \label{eq:complexstep}
\dot{\lambda} = -\frac{\text{Img}[H(\tilde{\textbf{x}})]}{\gamma_c},
\end{equation}
where $\tilde{\textbf{x}} = \textbf{x} + \gamma_c i$ denotes a complex version of the state vector in which the state associated with the costate has a complex component only, and $\gamma_c$ denotes the complex perturbation step size taken in the imaginary direction, $i$. Img() returns the imaginary part of the Hamiltonian. Very small complex step sizes can be taken, i.e., $\gamma_c = 1.0 \times 10^{-16}$. 

The accuracy of the above approach for obtaining the derivative of a function is, remarkably, shown to match the analytical partial derivatives with near machine precision \cite{martins2013review,lantoine2012hybrid} for a wide range of $\gamma_c$ values, in contrast to real-valued finite-difference approximations, where a ``Goldilocks'' step size region exists for obtaining valid partial derivative approximation. Note that the CX method evaluates the numerical value of the Hamiltonian. Obviously, \textit{numerical evaluation of the Hamiltonian is significantly easier than constructing the associated symbolic expression of the Jacobian. This is a key enabler when we deal with systems with complicated dynamics}.

The relation defined in Eq.~\eqref{eq:complexstep} is used at each time instant (say, in a standard Runge-Kutta method) to obtain the Jacobian of the Hamiltonian with respect to the state variables. There is a key point regarding the implementation of the CX-based method (or any other analogous numerical treatment). In order to use this method, recall that the Euler-Lagrange equation ($\dot{\lambda}_{x_i} = - \partial H/\partial x_i$, $i=1,\cdots,7$) sets the guideline for deriving the costate dynamics. In other words, the costate dynamics do not depend on the structure of the control that appears in the Hamiltonian. On the other hand, for constructing the numerical value of the Hamiltonian, information regarding both states and control is required. However, when one seeks to re-construct the \textit{costate dynamics numerically}, optimal control inputs (that also depend on states of the system through switching function and costates) \textit{should not} be modified at all. 

At this stage and in order to clarify a number of important points, we discuss the implementation of the method when the equations of motion are formulated using Cartesian coordinates. Comparisons to the analogous implementation are made when the set of MEEs is used. Further details and discussions of the resulting algebraic switching function expressions and form of the primer vector are given in \cite{taheri2016enhanced,junkins2018exploration}.

In fuel-optimal problems, if Cartesian coordinates are used to represent the dynamics, the unit vector, $\hat{\bm{\alpha}} = - \bm{\lambda}_{\textbf{v}}/\|\bm{\lambda}_{\textbf{v}} \|$, is a smooth function that solely depends on the costate vector associated with the velocity vector. This control is substituted in the Hamiltonian so that it is possible to determine the optimal policy for the thrust value (through switching function). The fact that optimal, $\hat{\bm{\alpha}}$ is not a function of the states, $\mathbf{x}$, would not create any problem when $\hat{\bm{\alpha}} = - \bm{\lambda}_{\textbf{v}}/\|\bm{\lambda}_{\textbf{v}} \|$ is substituted back into the equations of motion (which will eventually appear in the Hamiltonian). However, the throttle input is a function of the mass through the switching function 
\begin{equation}
S = \frac{\|\bm{\lambda}_{\textbf{v}} \|}{m}+\frac{\lambda_m}{c}.
\end{equation}

Recall that we need to evaluate the value of the Hamiltonian numerically, but the numerical evaluation of the costate dynamics should be consistent with the Euler-Lagrange equations. It is obvious that when $\dot{\lambda}_m$ is calculated numerically, the switching function will be affected through the mass in the denominator of the first term. Eventually, this will affect the numerical value of the Hamiltonian.

The situation is worse when the set of MEEs is used for dynamic modeling since both the optimal direction (through $\mathbb{B}(\textbf{x})$) and the switching function (through $\mathbb{B}(\textbf{x})$ and $m$) depends on the the state as
\begin{align}
\hat{\bm{\alpha}}^* = & -\frac{\mathbb{B}^{\top} \pmb{\lambda}}{\|\mathbb{B}^{\top} \pmb{\lambda} \|}, & S = & ~ \frac{\| \mathbb{B}^{\top} \pmb{\lambda} \|}{m}+\frac{\lambda_m}{c}.
\end{align}

Ultimately, when numerical approaches are to be used to re-construct the costate time derivative values numerically, \textit{the control has to be computed with the original state vector, i.e., with no complex component}. Then, to determine the gradient of $H$, the control inputs are treated as constants wherever they appear in $H$ expression while the other state-dependent terms will be affected by a perturbation in the imaginary component. 

For the problem under study, the Hamiltonian that is used in Eq.~\eqref{eq:complexstep} becomes
\begin{equation}
H(\tilde{\textbf{x}}) = \bm{\lambda}^{\top} \left \{ \textbf{A}(\tilde{\textbf{x}},t) +  \mathbb{B}(\tilde{\textbf{x}},t) \left [ \frac{2 \eta}{m~\boxed{c}} \boxed{\hat{\bm{\alpha}}} + \textbf{a}_{\text{sb}}(\tilde{\textbf{x}}) \right ] \right \} - \lambda_m \frac{2 \eta \boxed{P}}{(\boxed{c})^2}.
\end{equation}

The terms enclosed in boxes should not be updated when the CX-based method is used. Note also that for VIVT engines, $c$ is also a control input, which means that it depends on, $\textbf{x}$ as derived in Eqs.~\eqref{eq:VIVTcop} and \eqref{eq:VIVTcthree}. A salient feature of the numerical approach is that it is amenable to parallel implementation since each of the states of the system (in our case, seven states) can be perturbed (in the complex component) separately and independently; the above outlined procedure can be followed for each state. 

We have developed a parallel code to implement the CX-based method, but chose not to use it in this paper since the computations were found not demanding and this parallel tool is not necessary to demonstrate the main ideas. Even though we need to include the effects of gravitational perturbations due to secondary bodies, the introduction of these perturbations proves not to be problematic (for the considered test case) since we already have an excellent near-optimal solution to start from. The convergence when planetary gravity perturbations are considered is usually achieved in a small (single digit) number of iterations and takes only a few seconds for the example herein. The inclusion of gravitational perturbations due to secondary bodies, and when dynamics is modeled using the MEEs, with the CX-based method is another contribution of the paper. 

\section{Control Smoothing for VIVT Engines} \label{sec:CSC}
Recently, a generic smoothing approach has been developed, based on HTS \cite{taheriagenericsmoothing2018}. Its utility has been demonstrated for smoothing bang-off-bang and bang-bang control inputs on a number of problems including interplanetary and geocentric low-thrust trajectory optimization problems \cite{junkins2018exploration}, rest-to-rest attitude re-orientation of rigid-body and flexible satellites \cite{taheriagenericsmoothing2018,sandeep2018hybrid}, and low-thrust transfer from a geostationary transfer orbit to an L1 Halo orbit in the circular restricted circular three-body model of the Earth-Moon system \cite{taheri2018generic}. Its application for constructing a smooth transition between regular and singular control inputs has been readily demonstrated for the maximum-altitude Goddard rocket problem \cite{taheri2018doublesmoothing}. It is also used for constructing the \textit{optimal switching surfaces} that provide insights to optimal trajectories through which Taheri and Junkins \cite{taheri2018how} have answered the Edelbaum's famous question: ``How many impulses?'' \cite{edelbaum1967many}.

In these settings, the HTS method is used to replace the instantaneous sharp changes in control with a controllably sharp, rigorously continuous smooth approximation, where one smoothness parameter controls the sharpness of the switches. In this case, the homotopic smoothness parameter, $\rho_c$ (subscript `c' for constraint-related smoothing parameters) lies in the interval $\rho_c \in (0,\rho_{c,\text{max}}]$. The upper bound can initially be set to unity (or higher values) to maximize smoothness and to promote convergence; then a homotopic sweep of the parameter, $\rho_c$ toward zero is used to solve subsequent solutions. However, the solution is smoothly approaching the instantaneous switch limit as $\rho_c \rightarrow 0$. 

The application of the HTS method has been found to be quite an effective approach since it enlarges the domain of convergence of the resulting TPBVPs such that, in a large number of cases and, with a moderate number of random sets of initial costate guesses, convergence is usually achievable \cite{junkins2018exploration}. Note, when modern adaptive step-size integrators are employed, the continuous approximation of switches eliminates the necessity for high-precision root-solving for switch times, which is a significant advantage. For discontinuous controls with state- and costate-dependent switch times, note that conventional step-by-step integrators will not yield a valid integration step if a force discontinuity is contained interior to an integration step. Therefore, without smoothing, conventional integrators must be augmented with a logic to locate precisely the switch times and to use these times to adjust integration step size values \cite{saghamanesh2018robust}. 

In the present discussion, we introduce another variant of the HTS method to construct a composite control by blending a set of distinct controls continuously. The principal guideline is to use a smooth equivalent representation to replace discontinuous ``if $\rightarrow$ then'' conditions as these conditional statements are the root cause of control non-smoothness. A signed measure of distance from satisfaction of the ``if $\rightarrow$ then'' discontinuity is used as the argument of a smoothed activation function. In order to use the HTS method, we represent each of the ``if $\rightarrow$ then'' conditions as a continuous approximation of the constraint that triggers an activation function. 

Qualitatively, one can interpret the activation function to be a continuous sharp step function, which ultimately approaches unity when one (or a set of) constraint(s) condition is (are) satisfied or 0 when the constraint condition is violated. Another key point is that there could be a number of constraints for which, when all are nearly satisfied, a particular control ``dominates'' the other controls. Therefore, multiplicative incorporation of activation functions (each of which being weighted based on its corresponding constraint violation) allows us to smoothly approximate the switches. Eventually, a summation over all possible controls results in a composite continuous control that smoothly transitions centered on time points where, otherwise, instantaneous switches would occur.

For instance, let us consider the specific impulse relations defined in Eq.~\eqref{eq:VIVTcthree}. There are three possibilities for $c^*$ depending on the value of $c_{\text{op}}$, i.e., $c^* \in \{c_{\text{min}},c_{\text{op}}, c_{\text{max}}\}$. In fact, $c_{\text{min}}$ and $c_{\text{max}}$ are the limits on the value of specific impulse, which become active if the control violates these limits. A smooth composite representation of $c^*$ can be formed as
\begin{equation} \label{eq:VIVTcompositec}
c^*(\rho_\text{c}) = \zeta_{c_{\text{min}}} c_{\text{min}} + \zeta_{c_{\text{op}}} c_{\text{op}} + \zeta_{c_{\text{max}}} c_{\text{max}},
\end{equation}
where $\zeta_{c_{\text{min}}}$, $\zeta_{c_{\text{op}}}$, and $\zeta_{c_{\text{max}}}$ denote the activation functions corresponding to the lower, intermediate, and upper control inputs, respectively. The activation function for the lower control bound is defined as
\begin{align}
   \zeta_{c_{\text{min}}} = & \frac{1}{2} \left [ 1- \text{tanh} \left ( \frac{g_{c_{\text{min}},1}}{\rho_\text{c}} \right ) \right ], & g_{c_{\text{min}},1} = &~c_{\text{op}}-c_{\text{min}} \leq 0,
\end{align}
where one constraint, $g_{c_{\text{min}},1}$, has to be taken into consideration. Here, $\rho_\text{c}$ is the smoothing factor corresponding to the sharpness of the constraint-related step functions. The activation function for the intermediate control takes a multiplicative form, $\zeta_{c_{\text{op}}}$, which can be written as 
\begin{equation}
    \zeta_{c_{\text{op}}} = \zeta_{c_{\text{op}},1} \zeta_{c_{\text{op}},2},
\end{equation}
where two constraints are blended and their activation functions can be defined as
\begin{align}
\zeta_{c_{\text{op}},1} = & \frac{1}{2} \left [ 1- \text{tanh} \left ( \frac{g_{c_{\text{op}},1}}{\rho_\text{c}} \right ) \right ], & g_{c_{\text{op}},1} = &~c_{\text{min}}-c_{\text{op}} \leq 0,\\
\zeta_{c_{\text{op}},2} = & \frac{1}{2} \left [ 1- \text{tanh} \left ( \frac{g_{c_{\text{op}},2}}{\rho_\text{c}} \right ) \right ], & g_{c_{\text{op}},2} = &~c_{\text{op}}-c_{\text{max}} \leq 0. 
\end{align}

In a similar fashion, the activation function for the upper control, $\zeta_{c_{\text{max}}}$, can be written as
\begin{align} \label{eq:zetamax}
\zeta_{c_{\text{max}}} = & \frac{1}{2} \left [ 1- \text{tanh} \left ( \frac{g_{c_{\text{max}},1}}{\rho_\text{c}} \right ) \right ], & g_{c_{\text{max}},1} = &~c_{\text{max}}-c_{\text{op}} \leq 0.
\end{align}

It is straightforward to verify that Eq.~\eqref{eq:VIVTcompositec} is indeed a smooth representation of the discontinuous $c^*$ given in Eq.~\eqref{eq:VIVTcthree}. For instance, when $g_{c_{\text{min}},1} \leq 0$, then we know that $c=c_{\text{min}}$ is the solution based on Eq.~\eqref{eq:VIVTcthree}. According to the definition, $\zeta_{c_{\text{min}}} \rightarrow 1$ as $c_{\text{op}}$ gets smaller than $c_{\text{min}}$. At the same time, $\zeta_{c_{\text{op}}} \rightarrow 0$ and $\zeta_{c_{\text{max}}} \rightarrow 0$. Also, the smaller the value of $\rho_c$ is, the faster $\zeta_{c_{\text{min}}}$ approaches 1. The other controls become dominant if/when their associated constraint(s) is (are) triggered.
\begin{figure}[htbp!]
\centering
\includegraphics[width=4.0in]{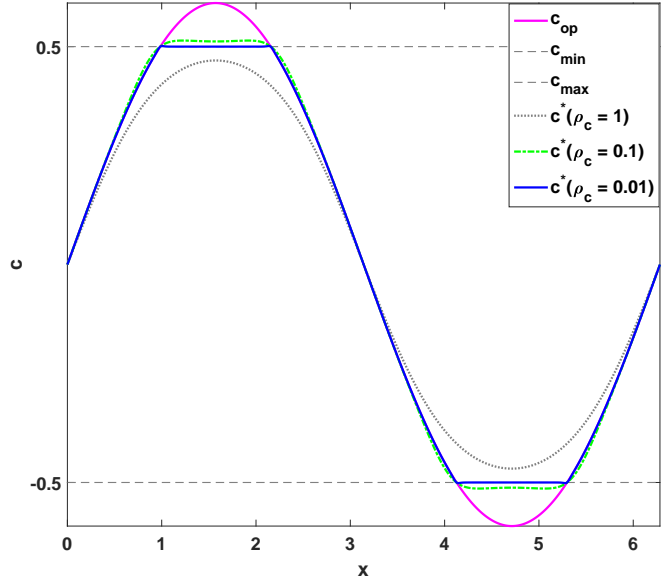}
\caption{Typical profiles of $c^*$ for different $\rho_c \in \{1.0, 0.1, 0.01\}$ values.}
\label{fig:TestcompositeC}
\end{figure}

A test problem is considered for demonstration. Let $x = [0,2\pi]$, it is assumed that $c_{\text{op}} = 0.6 \sin(x)$, $c_{\text{min}} = -0.5$ and $c_{\text{max}} = 0.5$. Figure \ref{fig:TestcompositeC} depicts the profiles of composite control, $c^*$ (using Eqs.~\eqref{eq:VIVTcompositec}-\eqref{eq:zetamax}) for different values of the smoothing parameter.


The other control input that we seek to smooth is the power input $P^*$ defined in Eq.~\eqref{eq:VIVTpowercondition}, which is usually idealized as having a bang-off-bang profile, depending on the value of $S$ defined in Eq.~\eqref{eq:sfm-p1}. Note that the bang-off-bang structure of power is one of the optimality criteria. In other words, even though the solar arrays are producing power (and there could be enough power to send to the engine), it is the power switching function defined in Eq.~\eqref{eq:VIVTpowercondition} that sets the criterion to switch the engine on or off. 

It is straightforward  to achieve a smooth approximation of such control types already illustrated in \cite{taheri2018generic} as
\begin{align}
P^*(S,\rho_{\text{b}}) = & \frac{1}{2} \left [ 1+ \text{tanh}\left ( \frac{S}{\rho_{\text{b}}} \right ) \right ] P^*_{\text{max}},\\ T^*(S,\rho_{\text{b}}) = & \frac{1}{2} \left [ 1+ \text{tanh}\left ( \frac{S}{\rho_{\text{b}}} \right ) \right ] T^*_{\text{max}},
\end{align}
where $\rho_{\text{b}}$ is the smoothing factor for bang-bang type control. The possibility of having singular-type power arcs is not considered in this work as they rarely occur in optimal space flights.

In summary, in the TPBVP associated with the fuel-optimal trajectory of a VIVT thruster, the following smooth approximations are implemented in the RHS of the set of state-costate dynamic equations defined in Eq.~\eqref{eq:F} with
\begin{align} \label{eq:power_and_c_op}
P^*(S) & \approx P^*(S,\rho_{\text{b}}), & c^* & \approx c^*(\rho_{\text{c}}).
\end{align}

While we have introduced the HTS smoothing heuristically above, it will be evident that the optimal instantaneous switches associated with Pontryagin's necessary conditions for optimality can be smoothly approximated by making the values of $\rho_b$ and $\rho_c$ sufficiently small. Note that there are two different smoothing parameters and their values can be different depending on the type of the problem. More details regarding the difference and judicious selection of these parameters are given in the Section \ref{sec:results}. 

\subsection{Generalization of Composite Smooth Controls} \label{sec:generalcomposite}
Let $m \in \mathbb{N}$ denote the number of scalar control inputs of a system. Let $n_{b,k}$ denote the number of building-block controls that have to be merged smoothly to fully define the $k$-th control time history, $1 \leq k \leq m$, and let $u_{b,i}$ denote the $i$-th building block control. 
The composite smooth representation of the $k$-th control input can be written as
\begin{equation} \label{eq:generalizedform}
u_{k} = \sum_{i=1}^{n_{b,k}} \left [~\prod_{j=1}^{n_{c,k,i}} \zeta_{i,j} \right ] u_{b,i}, ~~~~~ k = 1, \cdots, m,
\end{equation}
where
\begin{align} \label{eq:conact_smooth}
\zeta_{i,j} & = \frac{1}{2} \left [ 1- \text{tanh}\left ( \frac{g_{i,j}}{\rho_{\text{c}}} \right ) \right ],
\end{align} 
and $n_{c,k,i}$ is the number of constraints on the $i$-th building block control, which are handled by the smoothing functions. Recall that the actual activation functions, $\zeta_{i,j}$, have the following forms
\begin{align} \label{eq:conact}
    \zeta_{i,j} & =\begin{cases}
 1, & \text{if} ~ g_{i,j} \leq 0,\\
 0, & \text{otherwise}.
\end{cases}
\end{align}
Equation.~\eqref{eq:conact_smooth}, which is substituted in Eq.~\eqref{eq:generalizedform}, represents a smooth approximation of the actual activation function defined in Eq.~\eqref{eq:conact}.

As an example, the case given in Eq.~\eqref{eq:VIVTcompositec} is reviewed. Here, $m = 2$ since smooth approximations of both $P^*$ and $c^*$ are sought. Let $k=1$ correspond to $c^*$ for which $n_{b,1} = 3$, as there are three building-block controls, i.e., $c^* \in \{c_{\text{min}},c_{\text{op}}, c_{\text{max}}\}$. Accordingly, there are different number of constraints for each control, i.e., $n_{c,1,1} = 1$, $n_{c,1,2} = 2$ and $n_{c,1,3} =1$. Using Eq.~\eqref{eq:generalizedform}, we can write 
\begin{equation} \label{eq:cstarcomposite2}
c^* = u_1 = \zeta_{1,1} u_{b,1} + \zeta_{2,1} \zeta_{2,2} u_{b,2} + \zeta_{3,1} u_{b,3}.
\end{equation}

By comparing Eq.~\eqref{eq:cstarcomposite2} with Eq.~\eqref{eq:VIVTcompositec}, it is easy to verify that 
\begin{align}
    \begin{cases}
    c_{\text{min}} = u_{b,1},\\
    c_{\text{op}} = u_{b,2},\\
    c_{\text{max}} = u_{b,3},
    \end{cases}  \text{and}~~\begin{cases}
    \zeta_{c_{\text{min}}} = \zeta_{1,1},\\
    \zeta_{c_{\text{op}}} = \zeta_{2,1} \zeta_{2,2},\\
    \zeta_{c_{\text{max}}} = \zeta_{3,1}.
    \end{cases}
\end{align}
\subsection{Consideration of No-Power Condition}
It is also possible to handle situations in which there is not sufficient power to energize the PPU. For solar power-driven electric propulsion, this can happen at farther distances from the Sun. Such a switch possibility can be characterized by the smooth approximation
\begin{equation} \label{eq:nopow}
\zeta_{\text{no-power}} = \frac{1}{2} \left [ 1+ \text{tanh} \left ( \frac{P_{\text{ava}}-P_{\text{ava,min}}}{\rho_{\text{b}}} \right ) \right ],
\end{equation}
where $P_{\text{ava,min}}$ is a user-defined minimum power limit. We have used the same smoothing parameter used for smoothing bang-bang power behavior, $\rho_b$ to avoid introducing another smoothing parameter. This constraint acts upon the power level and it is assumed that the same value of smoothing parameter can be used for both conditions. 

The activation function defined in Eq.~\eqref{eq:nopow} should be multiplied by the power approximation defined in Eq.~\eqref{eq:power_and_c_op} with
\begin{equation}
P^*(S) \approx \zeta_{\text{no-power}} P^*(S,\rho_{\text{b}}).
\end{equation}

\subsection{Imposing No-Thrust Time Interval and Thrust Vector Direction Constraints}
It is also possible to 1) enforce a no-thrust time interval constraint, and 2) impose a certain condition on a direction of interest for the thrust vector. Let $t_l$ and $t_u$ denote two time instants such that $t_0 < t_l <t_u <t_f$. It is then possible to define an activation function that is triggered only when $t \in [t_l, t_u]$ as
\begin{equation}
\zeta_{\text{zero-thrust}} = \zeta_{t_l} \zeta_{t_u},
\end{equation}
where
\begin{align}
\zeta_{t_l} = & \frac{1}{2} \left [ 1+ \text{tanh} \left ( \frac{t-t_l}{\rho_\text{c}} \right ) \right ], & \zeta_{t_u} = & 1-\frac{1}{2} \left [ 1+ \text{tanh} \left ( \frac{t-t_u}{\rho_\text{c}} \right ) \right ].
\end{align}

Eventually, the function $P^*(S)$ can be replaced by
\begin{equation}
P^*(S) \approx (1-\zeta_{\text{zero-thrust}})\zeta_{\text{no-power}} P^*(S,\rho_{\text{b}}).
\end{equation}

The capability to impose no-thrust time interval is helpful if one is interested to perform missed-thrust analyses \cite{laipert2015automated} and also to accommodate events such as time windows when sensitive science measurements will be made. Note that one is able to use any other type of state-dependent constraint to enforce a particular pattern on thrust. For instance, it is possible to handle shadow constraints by using a particular angle-like variable. Alternatively, once can incorporate shadow events by a multiplicative smooth activation function over a controllably small window centered on an event measured by a Sun-detection indicator. 

Similarly, during a given time period, $t\in[t_l,t_u]$, it is possible to constrain a particular specified (time-varying or fixed) direction of interest, $\hat{\bm{\alpha}}_{\text{desired}}$, for the thrust vector by replacing $\hat{\bm{\alpha}}^*$ in Eq.~\eqref{eq:primervector} with
\begin{equation} \label{eq:smooththrustdirection}
\hat{\bm{\alpha}}^* = (1-\zeta_{\text{time}}) \frac{\textbf{p}}{\|\textbf{p} \|} + \zeta_{\text{time}} \hat{\bm{\alpha}}_{\text{desired}},
\end{equation}
where 
\begin{equation}
    \zeta_{\text{time}} = \zeta_{t_l} \zeta_{t_u}.
\end{equation}

During the time interval, $t\in[t_l,t_u]$, the direction of the thrust vector is re-oriented away from the unit vector along the primer vector, $\textbf{p}$, toward $\hat{\bm{\alpha}}_{\text{desired}}$. 

\subsection{Connection Between Smooth and Non-smooth Necessary Conditions}
The particular type of homotopy that we have employed in this framework seeks to establish a smooth representation neighboring all non-smooth components of the resulting boundary-value problems. Specifically, real-valued inputs including the direction of thrust vector and power to the VIVT engine or power input to each engine and the number of active engines in a multi-mode cluster of engines (in part 2) have been replaced by smooth representations. Through choice of the smoothing parameter, approximations can be brought arbitrarily close to instantaneous (discontinuous) switches between a finite number of modes. Remarkably, this smooth embedding solution can be found through affordable computational process. 

A one- or multiple-parameter family of neighboring smooth OCPs are constructed (denoted by superscript `s') such that as the set of homotopy parameters (in our problem $\rho_b$ and $\rho_c$) approach a limiting value (0 or 1 depending on the formulation), the set of optimality conditions of the original MPBVP is recovered and satisfied. The smoothed problems are found to be much more computationally attractive than the formal optimal solution, which require the discontinuous multi-mode switching. 

We construct the smoothed version of the optimality conditions to ensure the following mapping (invariant embedding):
\begin{align}
    \begin{cases}
    \dot{\textbf{x}}^s = \dot{\textbf{x}}^s(\textbf{z}^s,\mathbf{U}^s,t)\\
    \dot{m}^s = \dot{m}^s(\textbf{z}^s,\mathbf{U}^s,t),\\
    \dot{\bm{\lambda}}^s= \dot{\bm{\lambda}}^s(\textbf{z}^s,\mathbf{U}^s,t),\\
    \dot{\lambda}_m^s = \dot{\lambda}_m^s(\textbf{z}^s,\mathbf{U}^s,t),
    \end{cases} \xrightarrow[\bm{\rho} = \bm{\rho}_{\text{max}} \rightarrow \mathbf{0}]{\mathbf{U}^s \rightarrow \mathbf{U}^* }
    \begin{cases}
    \dot{\textbf{x}}^* = \dot{\textbf{x}}^*(\textbf{z}^*,\mathbf{U}^*,t),\\
    \dot{m}^* = \dot{m}^*(\textbf{z}^*,\mathbf{U}^*,t),\\
    \dot{\bm{\lambda}}^* = \dot{\bm{\lambda}}^*(\textbf{z}^*,\mathbf{U}^*,t),\\
    \dot{\lambda}_m^* = \dot{\lambda}_m^*(\textbf{z}^*,\mathbf{U}^*,t),
    \end{cases}
\end{align}
where $\bm{\rho} = [\rho_b,\rho_c]^{\top}$ denote the vector of smoothing parameters. Smooth and optimal control vectors (that consist of all control variables) are given as 
\begin{align*}
\mathbf{U}^s \in & \left [ \bm{\alpha}^s(\bm{\rho}), P^*(\rho_b ), c^*(\rho_c) \right ],\\
\mathbf{U}^* \in  & \left [ \bm{\alpha}^*(\bm{\rho}_{\text{min}}), P^*(\rho_b = \rho_{\text{b,min}}), c^*(\rho_c = \rho_{\text{c,min}}) \right ].
\end{align*}

\indent The maximum values for the smoothing parameters, $\rho_{\text{b,max}}$ and $\rho_{\text{c,max}}$, are problem dependent, and their lowest values ($\rho_{\text{b,min}}$ and $\rho_{\text{c,min}}$) would have to be precisely zero to recover the original solution. In practice, as we sweep $\rho_b$ and $\rho_c$ toward zero we can monitor $J(\bm{\rho})$ (Eq.~\eqref{eq:costfunction}) and usually find $\bm{\rho}_{\text{min}}$ such that ``engineering optimality'' is achieved for non-zero $\bm{\rho}$ values. The control vector for the considered problem in this paper consists of only three components. However, it is straightforward to extend same idea to problems with additional control inputs and continuation parameters, which is treated in part 2.

On the other hand, we have found that the original structure of optimal control is usually recovered with negligible sacrifice on optimality if, $\bm{\rho}$ is not precisely zero (e.g., see \cite{taheri2016enhanced,taheri2018generic}). In fact, for a large category of problems, a smooth approximation of theoretically sharp control inputs are far more practical. A prominent example of such appealing features is the reorientation of flexible satellites. Sharp control inputs excites unmodeled dynamics and the use of ``smoothed versions'' of instantaneous switch control profiles is ideal to perform a realistic maneuver \cite{junkins1993introduction}. 
\section{Numerical Results} \label{sec:results}
To demonstrate the utility of the proposed framework, a low-thrust trajectory from Earth to a potentially hazardous asteroid (Dionysus) is considered \cite{taheri2016enhanced}. Due to the large change between the inclination, eccentricity, and semi-major axis orbital elements of the Earth and those of the asteroid Dionysus, low-thrust trajectories would usually consist of several revolutions around the Sun and, depending on mass, $I_{\text{sp}}$ and propulsion system design, may take up to several years. 

Consequently, it is a challenging problem for the conventional indirect optimization methods. Furthermore, the effects of solar array degradation and perturbations due to secondary bodies can become significant over such prolonged multi-year, multi-revolution maneuvers. 

The terminal position and velocity vectors are taken from \cite{taheri2016enhanced} with one slight modification. Since we are taking into account the perturbing acceleration due to all planets, the initial position vector cannot match that of the Earth. Therefore, it is assumed that the spacecraft is on the boundary of the SOI of the Earth with a positive along-the-track offset of one Earth's SOI radius. Therefore, the initial position and velocity vectors of the spacecraft at $t = t_0$ are
\begin{align}
\textbf{r}_0 & = \begin{bmatrix}
-4561588.65006029\\
147076954.664376\\
-2259.94592436179
\end{bmatrix} \text{km}, & \textbf{v}_0 = \begin{bmatrix}
-30.2650979882182\\
-0.848685467901138\\
5.05303606281563 \times 10^{-5}
\end{bmatrix} \text{km/s}.
\end{align}
The final position and velocity vectors are
\begin{align}
\textbf{r}_f & = \begin{bmatrix}
-305026788.667814\\
307051467.941918	\\
82899899.5682193
\end{bmatrix} \text{km}, & \textbf{v}_f = \begin{bmatrix}
-4.23872656978066\\
-13.436307899221\\
0.565362569286115
\end{bmatrix} \text{km/s}.
\end{align}

The time of flight is fixed at $t_f - t_0 = 3543$ days. The initial mass of the spacecraft is $m_0 = 4000$ kg, its VIVT thruster has an assumed constant efficiency of $\eta = 0.65$. It is also assumed that the engine design is such that $I_{\text{sp}}$ can vary between 3000 and 6000 seconds. The beginning-of-life power is set to $P_{0,\text{BOL}} = 10$ kW and $\sigma = 2$\% solar array degradation per year is considered. It is assumed that 400 watts of power are needed to energize the PPU and operate various spacecraft sub-systems during the entire time of flight. 

Planetary perturbations in the modeling include the disturbing acceleration due to all of the planets of the solar system modeled through Eq.~\eqref{eq:secondaryacc}, from the innermost planet Mercury to the outermost planet Neptune, $N_{\text{sb}} = 8$. The numerical CX-based method is used to evaluate the contribution of planetary perturbations into costate dynamics. 

To quantify the impact of various models, and given the flexibility of the tool developed in this work, three cases are considered and are listed in Table \ref{tab:ED_VIVT_Cases}. The difference in these cases is due to the inclusion of degradation in the efficiency of the solar arrays and planetary perturbation models summarized as follows
\begin{itemize}
    \item \textbf{Case 1}: two-body gravitational model without  consideration of variation of power due to change in distance and degradation of the solar arrays, and no inclusion of planetary perturbations,
    \item \textbf{Case 2}: two-body gravitational model with consideration of variation of power due to change in distance and degradation of the solar arrays, and no inclusion of planetary perturbations,
    \item \textbf{Case 3}: two-body gravitational model with consideration of variation of power due to change in distance and degradation of the solar arrays, and with inclusion of all planetary perturbations.
\end{itemize}

Compared to other cases, case 1 results in the maximum value of the final mass due to the fact that solar arrays' efficiency degradation due to time and planetary perturbations (albeit not very  significant) are entirely ignored. We can see that the solar power system degradation of 2\% per year ``costs'' about 38 kg of payload in this case. The difference between the final masses in cases 2 and 3 indicates that planetary gravitational perturbations are of secondary significance as pointed out in \cite{rayman2002design}. While the planetary perturbations have a small impact on the final mass (propellant consumption), they nonetheless result in trajectory deviations that must be included for a real mission design. 
\begin{table}[htbp!] 
\begin{center} 
		\caption{Considered cases for the VIVT thruster for the Earth-to-Dionysus problem; $\rho_{\text{b}} = \rho_{\text{c}}= 1.0 \times 10^{-5}$. } \label{tab:ED_VIVT_Cases}
		{\small
		\begin{tabular}{c c c c c c}
        \hline
        \hline
         Case  & Two-body & Power  & Degradation  & Planetary    & $m_f$\\
          \#       &          & Model  & Model        & Perturbation & (kg) \\    
                 &  $\mu_{\text{sun}}\textbf{r}/r^3 $        &  $\phi(r) = 1/r^2$  &  $\psi(t)$        & $\textbf{a}_{\text{sb}}$ & \\
         \hline          
          1      & Yes      & Yes         & No                & No  &  2848.1426  \\
          2      & Yes      & Yes         & Yes                & No  & 2786.2428 \\
          3      & Yes      & Yes         & Yes                & Yes & 2786.2383 \\
        \hline
        \hline
        \end{tabular}
		}
	\end{center}
\end{table}

\begin{figure}[htbp!]
\centering
\includegraphics[width=3.5in]{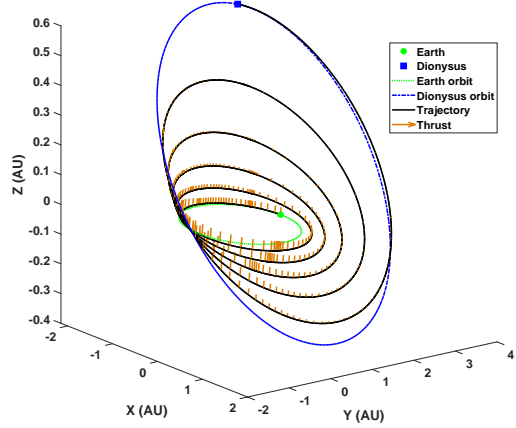}
\caption{Three-dimensional view of the Earth-to-Dionysus trajectory for case 3.}
\label{fig:VIVT_ED_TrajCase3}
\end{figure}

Figure \ref{fig:VIVT_ED_TrajCase3} depicts the location of the Earth in its orbit at the time of departure (December 23, 2012), low-thrust trajectory, and location of the asteroid in its orbit at the end of flight, all in the heliocentric J2000 frame of reference. The optimal solution corresponds to making five revolutions around the Sun. In this respect, there is no change to the ``optimal'' number of revolutions when a different type of thruster is used \cite{taheri2016enhanced}. In fact, the amount of propulsive force generated with the VIVT thruster is approximately on the same scale of the constant-$I_{\text{sp}}$ engine with throttling capability that is investigated in \cite{taheri2016enhanced}; hence, we observe no change to the optimal number of revolutions (compared to the optimal solution reported in \cite{taheri2016enhanced}). This is in direct relation to the general fact that number of revolutions is a strong function of the overall propulsive acceleration. 


\begin{figure}[htbp!]
\centering
\includegraphics[width=4.0in]{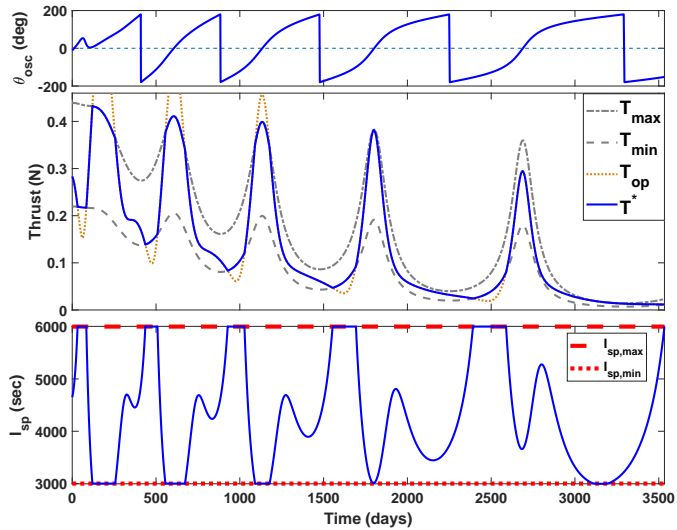}
\caption{Time history of osculating true anomaly, thrust, and specific impulse for case 3.}
\label{fig:VIVT_ED_TA_Thrst_IspCase3}
\end{figure}

Figure \ref{fig:VIVT_ED_TA_Thrst_IspCase3} depicts the time history of the osculating true anomaly, thrust level, and the specific impulse. The plot that shows the time history of thrust contains further information. There are three additional curves in the plot associated with the bounds corresponding to different thrust levels, i.e., $T_{\text{max}} = 2 \eta P_{\text{ava}}/c_{\text{min}}$, $T_{\text{min}} = 2 \eta P_{\text{ava}}/c_{\text{max}}$, and $T_{\text{op}} = 2 \eta P_{\text{ava}}/c_{\text{op}}$ (see Eq.~\eqref{eq:avapower}). The main difference between these thrust values is due to the value of $I_{\text{sp}}$ in the denominator. In fact, $T_{\text{max}}$ and $T_{\text{min}}$ define the boundaries of the envelope of optimal thrust values, $T^*$, while the profile of $T_{\text{op}}$ provides not only the optimal intermediate thrust values, but also the ``optimal transition'' between the boundaries. Note that the profile of $T_{\text{op}}$ (denoted by the bolder line) crosses these boundaries; however, thrust values beyond the limits are never adopted due to the limit on the value of the specific impulse (due to the expression developed in Eq.~\eqref{eq:VIVTcompositec}). 

The osculating true anomaly is plotted to show the regions during which the maximum thrust value is used. Clearly, the thruster switches to high-thrust mode at the perihelion passages ($\theta_{\text{osc}} \approx 0$) of intermediate quasi-elliptical orbits. In a similar fashion, the specific impulse switches to its lowest value during the maximum thrust arcs. As the spacecraft gets farther from the Sun, the interval during which the specific impulse assumes intermediate values gets larger.  

\begin{figure}[htbp!]
\centering
\includegraphics[width=4.0in]{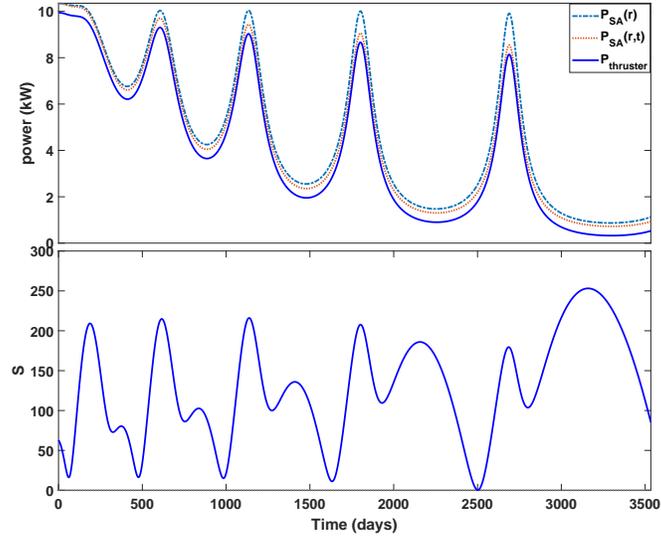}
\caption{Time histories of the solar arrays powers and the effects of time- and distance-dependent terms (upper plot) and the power switching function (lower plot) for case 3.}
\label{fig:VIVT_ED_PowersSwitchCase3}
\end{figure}

\begin{figure}[htbp!]
\centering
\includegraphics[width=4.0in]{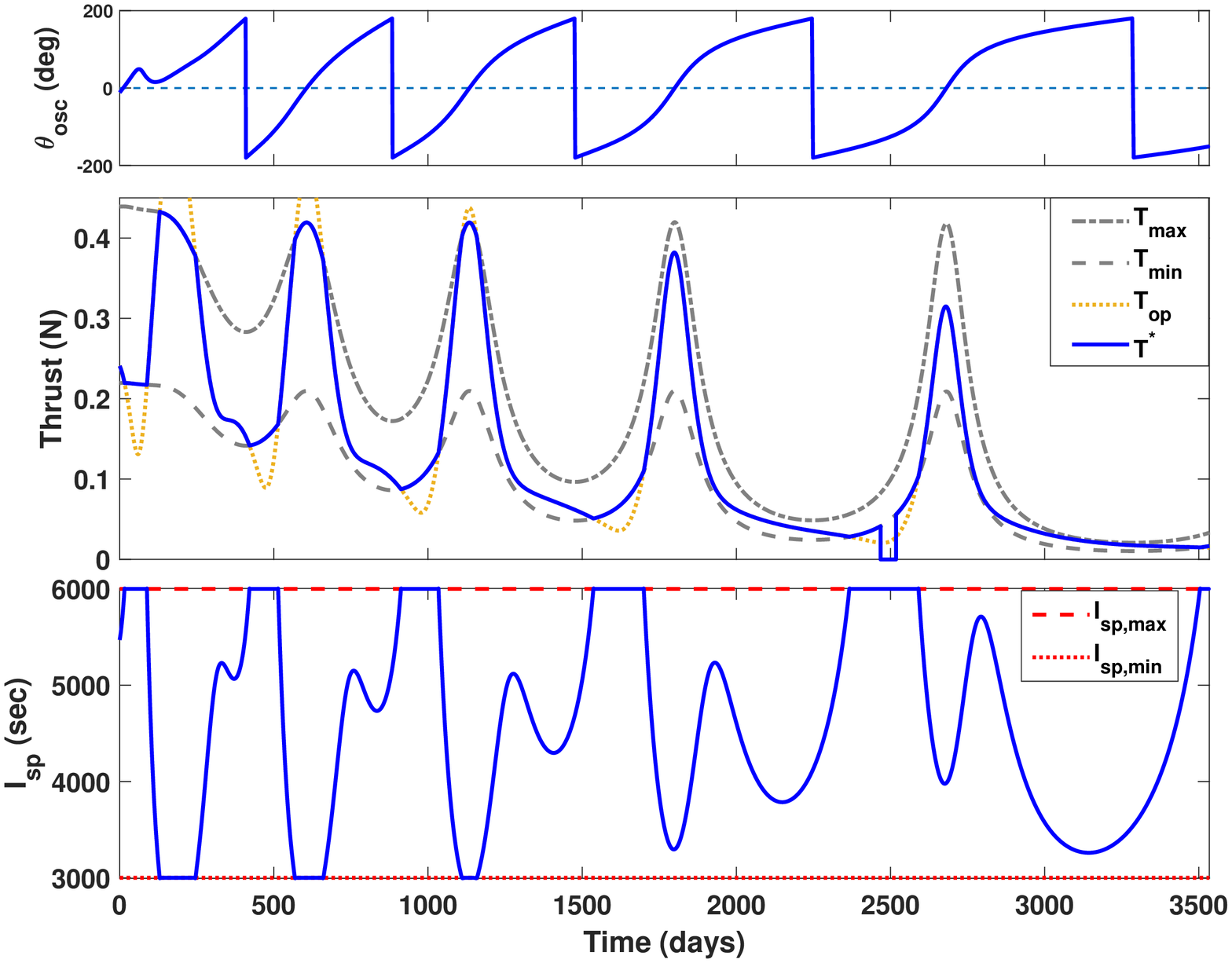}
\caption{Time history of osculating true anomaly, different thrusts, and $I_{\text{sp}}$ for case 1.}
\label{fig:VIVT_ED_TA_Thrst_IspCase1}
\end{figure}

Figure \ref{fig:VIVT_ED_PowersSwitchCase3} depicts the time history of power produced by the solar arrays and the power switching function. As the trajectory evolves, the degradation of solar arrays is more pronounced. The trend of the variation of power with time is expected as the distance to the Sun increases, with a consequent drop-off in solar intensity. The difference between the net power generated by the solar panel and the one sent to the thruster is always 400 watts. While not incorporated in the considered test cases, the CSC is capable of considering a maximum limit on the power produced by the solar arrays.
The power switching function remains entirely positive along the trajectory, which means that the thruster is operating continuously along the trajectory. The least value of the switching function occurs at $t = 2500$ days where it almost touches the $S = 0$ line. It is possible to have a no-thrust coast arc in the vicinity of this time with a slight variation of the system design parameters. Figure \ref{fig:VIVT_ED_TA_Thrst_IspCase1} shows that the no-thrust condition can occur for case 1 since the power switching function becomes negative during the no-thrust interval. A comprehensive study on the power switching function and its local features can be performed, similar to the one conducted in \cite{taheri2018how}.

\begin{figure}[htbp!]
\centering
\includegraphics[width=3.5in]{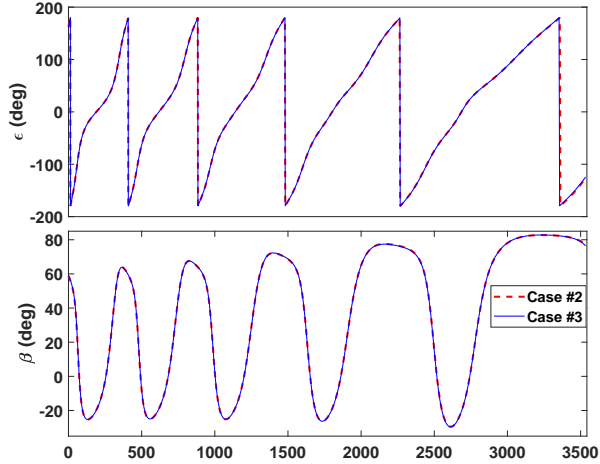}
\caption{The in-plane and out-of-plane thrust steering angles vs. time for cases 2 and 3.}
\label{fig:VIVT_ED_ThrustAnglessCase23}
\end{figure}

\begin{figure}[htbp!]
\centering
\includegraphics[width=3.5in]{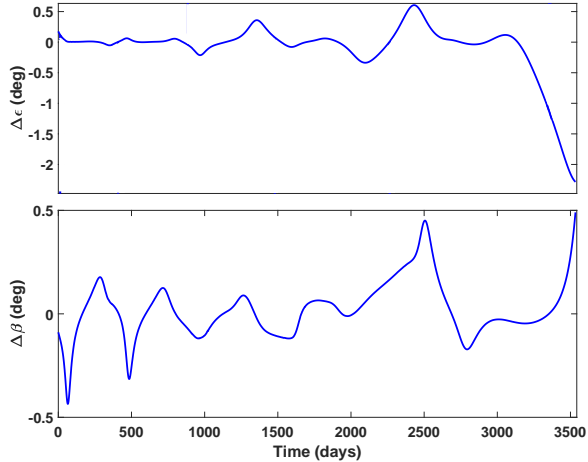}
\caption{Time history of the difference of the thrust steering angles between cases 2 and 3.}
\label{fig:VIVT_ED_ThrustAnglessCase23Dif}
\end{figure}

The results indicate that accounting for power sub-system dynamics in the model has a more significant impact on the final delivered mass than, in this case, planetary perturbations. However, the inclusion of perturbing accelerations due to secondary bodies will slightly alter the direction of the thrust vector along the trajectory. In order to demonstrate this point, we can parameterize the direction of the thrust vector (in the LVLH frame) by two angles $\epsilon$ and $\delta$. The former defines the in-plane angle between the projection of the thrust vector onto the $\hat{\textbf{u}}_{r}-\hat{\textbf{u}}_{t}$ plane of the LVLH frame, $\bm{\alpha}_{\text{proj}}$, and the $\hat{\textbf{u}}_{t}$ vector (measured positively clock-wise in the direction of $\hat{\textbf{u}}_{t}$ vector). The latter is the out-of-plane angle between the thrust vector $\bm{\alpha}_{\text{proj}}$ (measured positively in the direction of the specific angular momentum vector). Figure \ref{fig:VIVT_ED_ThrustAnglessCase23} depicts the time history of the in-plane and out-of-plane control steering angles for cases 2 and 3. Figure \ref{fig:VIVT_ED_ThrustAnglessCase23Dif} shows that there exist deviations as large as 2 degrees and 0.5 degrees in the in-plane and out-of-plane angles, respectively, solely due to the perturbations from the secondary bodies. 

\begin{figure}[htbp!]
\centering
\includegraphics[width=4.5in]{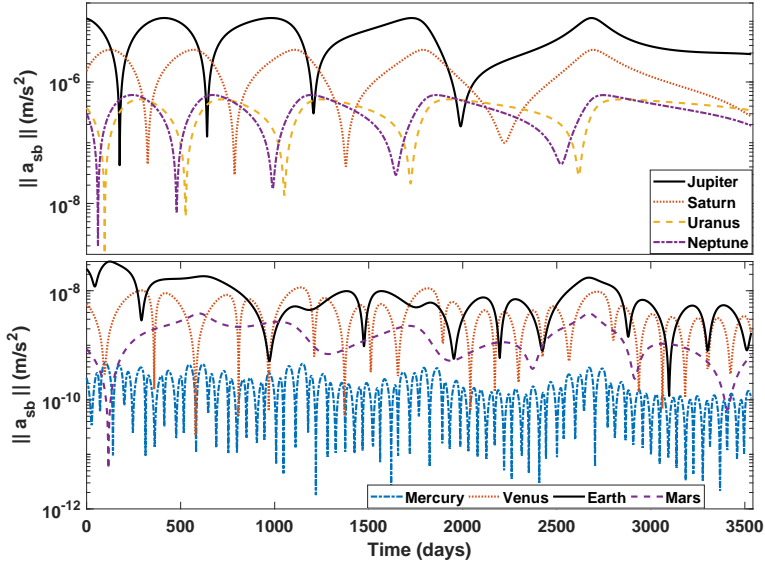}
\caption{The magnitude of the perturbing accelerations due to planets vs. time for case 3.}
\label{fig:VIVT_ED_SB_Normac}
\end{figure}

Figure \ref{fig:VIVT_ED_SB_Normac} shows the time history of the magnitude of the perturbing acceleration of the secondary bodies. The upper (lower) sub-plot shows the magnitude of the acceleration due to outer (inner) planets. As expected, the largest planets create greater magnitudes of disturbing accelerations, whereas among the inner planets, Earth creates the largest perturbation. Of course, these perturbations change drastically in the event that the trajectory makes a close encounter with any of the planets \cite{whiffen2006mystic}.

Obviously, the final mass is not a sole indication of the advantage of using a VIVT engine. The actual advantage of using any type of engine or a cluster of engines is determined through a multi-disciplinary design optimization, where the mass of the power and propulsion sub-systems, solar arrays, propellant and their contribution to the overall mass of the spacecraft are accounted for.

Our goal, in this paper, is to demonstrate the utility of the CSC framework to make the indirect optimization method for such problems amenable to numerical treatment. In addition, incorporation of realistic models for the power sub-system, efficiency degradation in solar arrays, and planetary perturbations models combined with the CX-based method provides us with significant flexibility to handle OCPs with complex state dynamics. The CSC framework enables us to solve OCPs using single-shooting schemes and MATLAB's \textit{fsolve} solver. We recognize that more sophisticated solvers and integrators and the use of compiled codes (using C++ or Fortran) accelerate these computations. In part 2 of this series of papers, application of the CSC framework is investigated for trajectory optimization of a spacecraft with multiple engines.

\section{Conclusion} \label{sec:conclusion}
In this work, a framework, called composite smooth control (CSC), is developed to deal with the problem of designing optimal control strategies for problems with multiple modes of operation. In particular, spacecraft equipped with a variable-$I_{\text{sp}}$ variable-thrust (VIVT) engine is considered for demonstration purposes. It is possible to modulate the value of the specific impulse to gain optimum efficiency, whereas the optimal switches between possible modes of operation are not known \textit{a priori}.  

Determination of the optimal transition times between the possible modes of operation of a VIVT engine is not a trivial task. The proposed CSC framework facilitates the numerical solution of the resulting multi-point boundary-value problems and allows us to treat them as two-point boundary-value problems with continuous control inputs, which is simpler to solve. 

For fuel-optimal problems, and for heliocentric phase of flight, the results indicate that accounting for power sub-system dynamics in the model has a more significant impact on the final delivered mass than, in this case, planetary perturbations. In the case of a VIVT engine, we have considered a complex multi-year multi-revolution interplanetary trajectory for which the thruster demonstrates a complex structure that consists of multiple modes of operation (i.e., values for specific impulse): maximum, minimum, and intermediate values. The profile of the specific impulse consists of 16 switches. In addition, implementation of the complex-based derivative approach is explained, where gravity perturbations due to all planets of the Solar System are taken into consideration, and when the set of modified equinoctial orbital elements are used for representing the dynamics.  

Moreover, the CSC framework has the capability to incorporate shadow- and time-triggered constraints in a smooth, continuous manner. These constraints will result in no-thrust arcs. The CSC framework enables us to impose time intervals during which the thrust vector has a possibly specific time variable direction. Thus, the methodology presented is promising for expanding the class of trajectory optimization problems to which indirect optimal control formulations can be successfully applied. 

\section{Acknowledgment}
Ehsan Taheri and John Junkins would like to thank our sponsors: AFOSR (Dr. Stacie Williams) and AFRL (Dr. Alok Das) for their support and collaborations under various contracts and grants. Ehsan Taheri and Ilya Kolmanovsky would also like to acknowledge the support of the National Science Foundation under the Award
Number CNS 1544844. 

\bibliographystyle{elsarticle-num}
\bibliography{References.bib}

\end{document}